\documentclass[12pt,leqno]{amsart}

\def\opn#1#2{\def#1{\mathop{\kern0pt\fam0#2}\nolimits}}

\opn\arg{arg} \opn\DD{\mathcal D} \opn\GCD{GCD}
\opn\Id{Id} \opn\im{Im} \opn\rad{rad}
\opn\sgdeg{sgdeg}  \opn\Frac{Frac} \opn\Spec{Spec}

\newtheorem{thm}{Theorem}

\newtheorem{lem}{Lemma}
\newtheorem{de}{Definition}
\newtheorem{nota}{Notations}
\newtheorem{rem}{Remark}

\newtheorem{rems}{Remarks}

\def\:{:\text{ }}
\def\CC{\mathbb C}

\def\fM{\mathfrak m}
\def\fP{\mathfrak P}

\def\fP{\mathfrak P}
\def\FF{\mathbb F}

\def\NN{\mathbb N}
\def\ZZ{\mathbb Z}

\begin{document}

\vspace*{3cm}
\begin{center}{\large\sc On the equivalence of the Jacobian, Dixmier and Poisson Conjectures in any characteristic}
\\[0.3cm]
 by \\[0.3cm]
{\bf Kossivi Adjamagbo \and Arno van den Essen}
\end{center}

\vspace*{1cm}

{\tiny {\bf Abstract:} Jacobian Conjecture in characteristic $p \geq 0$ \cite{keller} \cite{smale} \cite{adja1} means that for any positive integer $n$, any jacobian 1 endomorphism of the algebra of polynomials in $n$ indeterminates over a field of characteristic $p$ is an automorphism, provided it induces a field extension of degree not a multiple of $p$. Dixmier Conjecture in characteristic $p$ \cite{dixmier} means that, for any positive integer $n$, any endomorphism of the $n$-th Dirac quantum algebra over a field of characteristic $p$ \cite{dirac} \cite{littlewood} \cite{coutinho}, unjustly called Weyl algebra by J. Dixmier in \cite{dixmier}, i.e. the associative algebra over this field with $2n$ generators satisfying the normalized famous commutation relations of quantum mechanics, i.e. in other terms the algebra of ``formal'' differential operators in $n$ indeterminates with polynomials coefficients over this field \cite{adja0}, is an automorphism, provided its restriction to the center of this Dirac algebra induces a field extension of degree not a multiple of $p$, and the jacobian of this restriction is a non zero element of the field in the case where $p\leq n$. Poisson Conjecture in characteristic $p$ means that for any positive integer $n$, any endomorphism of the $n$-th canonical Poisson algebra over a field of characteristic $p$ , i.e. the algebra of polynomial in $2n$ indeterminates over this field endowed with its classical Poisson bracket, is an automorphism, provided it induces a field extension of degree not a multiple of $p$, and its jacobian of is a non zero element of the field in the case where $p\leq n$. Thanks to recent results on ring homomorphisms of Azumaya algebras \cite{ace} and to the following ones about endomorphisms of canonical Poisson algebras and Dirac quantum algebras, and about the reformulation in positive characteristic of these conjectures in characteristic zero on the model of \cite{adja2}, we prove the equivalence of these three conjectures in any characteristic, giving also by this way thanks to \cite{ace} a new proof of the equivalence of the complex version of the two first conjectures recently proved by Y. Tsuchimoto in a series of two papers \cite{tsuchimoto1} and \cite{tsuchimoto2}}

\vspace*{1cm}

\section{Introduction}

One of the simplest statements which is also one of the major mathematical problems for the new century, according to Steve Smale \cite{smale}, challenging the whole mathematical communauty for more than six decades, is the Jacobian Conjecture, more precisely the Jacobian Conjecture in characteristic zero \cite{keller}. It asserts that any jacobian 1 endomorphism of any algebra of polynomials in a finite number of indeterminates over a field of characteristic zero is an automorphism, see for instance \cite{bass} and \cite{arno}. Its generalization in any characteristic $p \geq 0$ introduced in  \cite{adja1} asserts that for any positive integer $n$, any jacobian 1 endomorphism of the algebra of polynomials in $n$ indeterminates over a field of characteristic $p$ is an automorphism, provided it induces a field extension of degree not a multiple of $p$. 

One the other hand, one of the most singular statements of the whole mathematical literature which is also challenging this community for almost four decades is Dixmier Conjecture, more precisely Dixmier Conjecture in characteristic zero. It asserts that any endomorphism of a Dirac quantum algebra over a field of characteristic zero \cite{dirac} \cite{littlewood} \cite{coutinho}, unjustly called Weyl algebra by J. Dixmier in \cite{dixmier}, i.e. the associative algebra over this field with generators $p_1,\ldots,p_n,q_1,\ldots,q_n$ satisfying the normalized famous commutation relations of quantum mechanics $[p_i,q_j] = \delta_{i,j}$ for each $i$ and $j$, i.e. in other terms an algebra of differential operators with polynomials coefficients over such a field, is an automorphism, see for instance \cite{dixmier}, \cite{bass} and \cite{arno}. Its generalization in any characteristic $p \geq 0$ that we are introducing asserts that that, for any positive integer $n$, any endomorphism of the $n$-th Dirac quantum algebra over a field of characteristic $p$ \cite{dirac} \cite{littlewood} \cite{coutinho},
called by mistake Weyl algebra by J. Dixmier in \cite{dixmier}, i.e. the associative algebra over this field with $2n$ generators satisfying the normalized famous commutation relations of quantum mechanics, i.e. in other terms the algebra of ``formal'' differential operators in $n$ indeterminates with polynomials coefficients over this domain  and Th. 2.12, is an automorphism, provided its restriction to the center of this Dirac quantum algebra induces a field extension of degree not a multiple of $p$, and the jacobian of this restriction is a non zero element of the field in the case where $p\leq n$. 

A similar singular statement is what we are introducing as Poisson Conjecture, more precisely as Poisson Conjecture in characteristic zero. It asserts that any endomorphism of a canonical Poisson algebra over a field of any characteristic zero, i.e. the algebra of polynomial in an even number of indeterminates over this field endowed with its classical Poisson bracket, is an automorphism. Its generalization in any characteristic $p \geq 0$ that we are introducing asserts that any endomorphism of a canonical Poisson algebra over a field of any characteristic $p$, i.e. the algebra of polynomial in an even number of indeterminates over this field endowed with its classical Poisson bracket, is an automorphism, provided it induces a field extension of degree not a multiple of $p$, and its jacobian of is a non zero element of the field in the case where $p\leq n$.

It is well known since the publication of \cite{bass} in 1982 that Dixmier Conjecture of index $n$ in characteristic zero implies the Jacobian one in dimension $n$. On the other hand, it is known for experts since the preprint \cite{arno2} from 1991 (see also \cite{arno}, p. 264) that the Jacobian Conjecture in dimension $n$ and in characteristic zero implies a weak form of Dixmier Conjecture of index $2n$ in characteristic zero, asserting that any endomorphism of a filtered Dirac quantum algebra of index $n$ over a field of characteristic zero, endowed with its fitration of ring of differential operators, is an automorphism. 

But the expected equivalence in any characteristic of the Jacobian conjecture in any dimension, the full Dixmier Conjecture of any index and Poisson Conjecture in any dimension is a kind of ``wedding of simplicity and singularity'', which is an interesting result on its own, independently from the proof of one of the ``united conjectures'', in conformity with Grothendieck's vision of mathematical research work, expressed in \cite{grothendieck}, Deuxi\`eme Partie, L'enterrement (I), and  which is a wonderful and rare lesson of methodology and wisdom : ``Ten things which are only guested, none of which (let us say Hodge Conjecture) carrying the conviction, but which light up and complement each other, like if they are working towards a same still mysterious harmony, find in this harmony the virtue of clarity. Even if all the ten will turn to be false, the work which led to this temporary clarity has not been fruitless, and the harmony of which it let us catch a glimpse and to which it gave us access for a while, is not an illusion, but a reality, inviting us to know it. Only by this work we have been able to be intimate with this reality, this hidden and perfect harmony. When we know that things are right to be what they are, that our vocation is to know them, not to dominate them, then the day when a mistake rise is a happy day, as much as the day when a proof learns us, beyond any doubt, that this thing that we were guessing is really the faithful and authentic expression of reality itself. In any case, such a discovery comes as a reward for a work, and can not occur without it. Even if it would come only after years of efforts, or even if we would never know the end of the story, reserved others after us, work is its own reward, rich at each moment with what the very moment reveal to us''.

So waiting patiently for the end of the history of the ``united conjectures'', the aim of the present paper is to expose the work which led us to the proof of their equivalence, mainly to deduce in any characteristic Dixmier Conjecture from Poisson Conjecture. This proof of the equivalence of the ``united conjectures'' also gives by this way thanks to \cite{ace} a new proof of the equivalence of the complex version of the two first conjectures recently proved by Y. Tsuchimoto in a series of two papers \cite{tsuchimoto1} and \cite{tsuchimoto2}.

However, the fact that he dont consider the crucial Azumaya property of Dirac quantum algebras over fields of positive characteristic which is well known since \cite{revoy} in 1973 by poeple interested in Dirac quantum algebras, and his apparent discovery in \cite{tsuchimoto2} of the well know properties of fields ultra-products which is well exposed for instance in \cite{eklof} ``for algebrists'' since 1977, added to a mistake in the last sentence of \cite{tsuchimoto2}(see the remarks following Theorem 3 below), dont make his proof as clear as desirable. Futhermore, it projects no light on Poisson Conjecture.

On the other hand, the paper \cite{bel2} of A. Belov and M. Kontsevich contains some unproved statements, one of which is clearly false and is fatal for the strategy of proof in this paper of the equivalence of Jacobian and Dixmier Conjectures in characteristic zero. This wrong statement, at the begining of section 4, claims that ``for any finitely generated domain $R$, we may assume that for any prime p the ring $R/pR$ is either zero or a domain''. However, the polynomial $X^4 + 1$ in one indeterminate $X$ over $\ZZ$ is irreducible in $\ZZ[X]$, while its canonical image in $(\ZZ/p\ZZ)[X]$ is reducible in this ring for any prime number $p$, as Alain Kraus drew our attention on it. So, $R = \ZZ[X]/(X^4 + 1)\ZZ[X]$ is a counter-example to this unproved and false claim. 

So, according to the eventful history of the Jacobian Conjecture, full with right unproved statements like ``Segre lemma'', finally proved by S. Abhyankar and T.-T. Moh (see for instance \cite{bass}, Faulty proofs) and Appelgate-Onishi-Nagata theorem, finally proved by S. Nagata (see for instance \cite{nagata}), and with wrong unproved ones (see for instance \cite{bass}, Faulty proofs), it was for us a duty toward history to propose the present explicit, and clear proof of the equivalence of the three indicated conjectures, not only in characteristic zero, but also in any characteristic. 

The expected clarity benefits greatly from the lights of basic and intimate properties of ring homomorphisms of Azumaya algebras (see \cite{ace}), of endomorphisms of canonical Poisson algebras (see section 2) and of endomorphisms of Dirac quantum algebras (see section 3) inspirated from the cited paper \cite{bel2}, \cite{bel1}, \cite{bel3} and \cite{tsuchimoto2}, and of the  reformulation in positive characteristic of the three conjectures in characteristic zero (see sections 4, 5 and 6), thanks to Los Theorem in Model Theory and ``Gabber bound'' for the degree of the inverse of automorphisms of polynomial, Poisson and Weyl algebras (see Theorem 2 below), on the model of the first explicit reformulation in positive characteristic of the Jacobian Conjecture ten years ago in \cite{adja2}, Theorem 3.9. , of which the present paper should be considered as the achievement. 

Before getting into the heart of the matter, we would like to express our deep gratitude to M. Kontsevich and A. Belov without who the present paper would probably never be writen, and to Charles-Michel Marle and Alain Kraus of University Paris 6 and Jean-Yves Charbonnel of University Paris 7, for fruitfull discussions during the preparation of this paper.

\section{Endomorphisms of canonical Poisson algebras over commutative rings}

\begin{nota} .
\begin{enumerate}
\item
Throughout this paper $R$ denotes a commutative ring with $1$, $n$ a positive integer, $X_1,\ldots, X_{2n}$ indeterminates over $R$, $X = (X_1,\ldots, X_{2n})$ $R^{[n]}$ the $R$-algebra $R[X_1,\ldots,X_n]$, $R^{[2 n]}$ the $R$-algebra $R[X_1,\ldots, X_{2n}]$.
\item
If $\phi$ is a non zero endomorphism of $R^{[n]}$ (resp. $R^{[2 n]}$), then $\deg_X(\phi)$ or $\deg(\phi)$ denotes the integer $max \{\deg_X(\phi(X_i)) | 1\leq i\leq n\}$ (resp. $max \{\deg_X(\phi(X_i)) | 1\leq i\leq 2n\}$).
\item
$P_n(R)$ denotes the $n$-th canonical Poisson algebra over $R$, i.e. the $R$-algebra $R^{[2n]}$ endowed with the canonical Poisson bracket  $\{, \}$ such that for any elements $f$ and $g$ of $P_n(R)$, we have : $$\{f,g\} = \sum_{i=1}^n (\frac{\partial f}{\partial X_i}\frac{\partial g}{\partial X_{i+n}}-\frac{\partial f}{\partial X_{i+n}}\frac{{\partial g}}{{\partial X_i}})$$ So, an endomorphism of $P_n(R)$ is an endomorphism $\phi$ of the $R$-algebra $R^{[2n]}$ such that $\{\phi(f),\phi(g)\} = \{f,g\}$ for any elements $f$ and $g$ of $P_n(R)$.
\item 
For any $F = (F_1,\ldots, F_n) \in (R^{[n]})^n$, $JF$ denotes the jacobian matrix of $F$, i.e. $(\frac{\partial F_i}{\partial X_j})_{1\leq i,j\leq n}$, where $i$ is the row index and $j$ the column index. For any endomorphism $\phi$ of $R$-algebra $R^{[n]}$, $J\phi$ denotes its jacobian matrix, i.e. the matrix $JF$, with $F = (\phi(X_1),\ldots, \phi(X_n))$.
\end{enumerate} 
\end{nota}

\begin{de} .

Let $B = (e_1,\ldots,e_n)$ be the canonical basis of the free $R^{[2 n]}$-module $E = (R^{[2 n])^n}$, $B^* = (e_1^*,\ldots,e_n^*)$ the dual basis of $B$.
\begin{enumerate}
\item
The canonical symplectic form on $E$ is the bilinear form : $$\omega = \sum_{i=0}^n {{e_i^*}\wedge {e_{i+n}^*}}$$ where, for $1\leq i\leq n$, ${e_i^*}\wedge {e_{i+n}^*}$ is the alternating bilinear form on $E$ such that for any integers $p$ and $q$ such that $1\leq p< q\leq n$, $({e_i^*}\wedge {e_{i+n}^*})(e_p,e_q) = 1$ if $i = p$ and $q = i+n$ and $0$ otherwise.
\item
The pull-back of $\omega$ by an endomorphism $L$ of the $R^{[2 n]}$-module $E$ is the bilinear form : $$L^* \omega = \omega \circ (L,L)$$
\item
Such an endomorphism $L$ is said to be symplectic if $L^* \omega = \omega$.
\item
An endomorphism $\phi$ of the $R$-algebra $R^{[2n]}$ is said to be symplectic if the endomorphism $L_\phi$ of the $R^{[2 n]}$-module $E$ the matrix of which in $B$ is the transposed of the jacobian matrix of $\phi$ is symplectic.
\item
If $R$ is an infinite domain, a polynomial map from $R^{2 n}$ to $R^{2 n}$ defined by an unique endomorphism $\phi$ of the $R$-algebra $R^{[2n]}$ is said to be symplectic if $\phi$ is symplectic.
 
\end{enumerate}
\end{de}

\begin{lem} .

For any endomorphism $\phi$ of the $R$-algebra $R^{[2n]}$, then the matrix in $B$ of the bilinear form $L_\phi^*\omega$ on $E$ is : $$(\{\phi(X_i),\phi(X_j)\})_{1\leq i,j\leq n}$$ 

\end{lem}

\begin{proof}

Let's put $F_i = \phi(X_i)$ for $1\leq i \leq n$ and let's consider integers $p$ and $q$ in $\{1,\dots,n\}^2$. Let's observe that : $$L_\phi(e_p) = \sum_{j=1}^{2n} \frac{\partial F_p}{\partial X_j}e_j$$ It follows that : $$L_\phi^*\omega(e_p,e_q) = \sum_{i=0}^n {(e_i^*\wedge e_{i+n}^*)( L_\phi(e_p),L_\phi(e_q))} = \{F_p,F_q\}$$

\end{proof}

\begin{lem} .

For any endomorphism $\phi$ of the $R$-algebra $R^{[2n]}$, then the following statements are equivalent :
\begin{enumerate}
\item
$\phi$ is symplectic.
\item
$\{\phi(X_i),\phi(X_j)\} = \{X_i,X_j\}$ for ${1\leq i,j\leq n}$.
\item
$\phi$ is an endomorphism of $P_n(R)$. 
\end{enumerate}
\end{lem}

\begin{proof} The equivalence of (1) and (2) follows from the previous lemma, knowing that two bilinear forms on $E$ are equal iff they have the same matrix in $B$. The implication $(1) \Rightarrow (2)$ is trivial and the inverse implication follows from the fact that the Poisson $\{,\}$ bracket on $R^{[2n]}$ is bilinear, antisymetric and satisfies Leibnitz's rule, i.e. for any $f$, $g$ and $h$ in $R^{[2n]}$, $\{f,gh\} = \{f,g\}h + \{f,h\}g$.

\end{proof}

\begin{thm} .

For any endomorphism $\phi$ of the $R$-algebra $R^{[2n]}$, the following statements are equivalent :
\begin{enumerate}
\item
$\phi$ is symplectic.
\item
$\phi$ is an endomorphism of $P_n(R)$ 
\end{enumerate}
If in addition one of these statements is true and $n!$ is inversible in $R$, then $\det(J\phi)=1$.
\end{thm}

\begin{proof} The equivalence of (1) and (2) follows from the previous lemma. So, let's assume $\phi$ symplectic and let's put $v = e_1^*\wedge \ldots \wedge e_{2n}^*$, $\omega^n = \omega_1\wedge \ldots \wedge \omega_n^*$, where each $\omega_i = \omega$. Since $v$ is the standard volume form on $E$, it is well-known that $\omega^n = n! (-1)^{(n(n-1))/2}v$ (see for instance \cite{godbillon}, Exemple 1.4, p. 123). So, since $n!$ is an invertible element of $R$, it follows that $\omega^n$ is a volume on $E$. Futhermore, since $\phi$ is symplectic, i.e. $L_\phi^*\omega = \omega$, it follows that $L_\phi^*\omega^n = \omega^n$. On the other hand, since $\omega^n$ is a volume form, it is well-known that $L_\phi^*\omega^n = (\det(L_\phi))\omega^n$ (see for instance \cite{godbillon}, Exec. (ii), p. 21). Since $\{\omega^n\}$ is a basis of the $R^{[2n]}$-module of all alternating $2n$-form on $E$, it follows that $\det(L_\phi) = \det(J\phi) =1$, as desired.

\end{proof}

\section{Endomorphisms of Dirac quantum algebras over commutative rings}

\begin{nota} .
\begin{enumerate}
\item
$A_n<R>$ denotes the $n$-th free algebra over $R$ with generators $Z_1,\ldots, Z_{2n}$, $\Im(A_n<R>)$ the bilateral ideal of $A_n<R>$generated by $[Z_i,Z_j] - \delta_{i+n,j}$ for $1\leq i, j\leq n$.
\item
$A_n(R)$ denotes the $n$-th Dirac quantum algebra over $R$, i.e. the $R$-algebra $$A_n<R>/\Im(A_n<R>) = R[Y_1,\ldots, Y_{2n}]$$ such that $Y_i = \rho(Z_i)$ for $1\leq i\leq n$, where $\rho$ is tha canonical map from $A_n<R>$ to $A_n<R>/\Im(A_n<R>)$, i.e. again the associative $R$-algebra with $2n$ generators $Y_1,\ldots, Y_{2n}$ and relations : $$For\  1\leq i, j\leq n, \  [Y_i,Y_j]=\delta_{i+n,j}$$ 
\item
If $a$ is a non zero element of $A_n(R)$, then $\deg(a)$ denotes the total degree of $a$ with respect to $Y_1,\ldots, Y_{2n}$, $Y = (Y_1,\ldots, Y_{2n})$ and if $\phi$ is a non zero endomorphism of $A_n(R)$, then $\deg_Y(\phi)$ or $\deg(\phi)$ denotes the integer $max \{\deg_Y(\phi(Y_i)) | 1\leq i\leq 2n\}$.
\item
If $R$ has a prime characteristic $p$, then it follows from Revoy Theorem 2 \cite{revoy} that $A_n(R)$ is an Azumaya algebra over its center denoted by $Z(A_n(R))$ and that this center is $R[Y_1^p,\ldots,Y_{2n}^p]$, which we identify with $R^{[2 n]}$ by choosing $X_i = Y_i^p$ for each $i$.
\item
If $R$ has a prime characteristic $p$, and if $\phi$ is an endomorphism of $A_n(R)$, then it follows from Revoy Theorem \cite{revoy} and from the Main Theorem of \cite{ace} that $\phi$ induces an endomorphism  on $R^{[2 n]}$ denoted $\phi_0$.
\end{enumerate}
\end{nota}

\begin{rems} . 
\begin{enumerate}
\item
A Dirac quantum algebra must not be confused with a ``Dirac algebra associated to a real vector space'' and generating a ``Clifford algebra'', as defined in Relativity theory (see for instance \cite{souriau}, p. 426-438).
\item
If $R$ has a prime characteristic $p$, and if $\phi$ is an endomorphism of $A_n(R)$, then it easy to observe that $\deg(\phi) = \deg(\phi_0)$ since, for $1\leq i\leq n$, we have : $$p\deg_X(\phi_0(X_i)) = \deg_Y(\phi_0(X_i)) = \deg_Y(\phi_0(Y_i^p)) =  $$ $$\deg_Y(\phi(Y_i^p)) = \deg_Y(\phi(Y_i)^p) = p\deg_Y(\phi(Y_i))$$
\item
The following theorem is a non commutative generalization of Gabber's degree bound theorem for the inverse of an automorphism of an algebra of polynomials over a field (see \cite{bass}, Cor. 1.4 of Theorem 1.5), completing the geometric generalization of this theorem for the inverse of an isomorphism of affine domains over a field (see \cite{adjawin}, Th. 3) :
\end{enumerate}
\end{rems}

\begin{thm}[degree bound theorem for automorphisms of Dirac quantum algebras] .

If $R$ is a commutative domain and $\phi$ an automorphism of the Dirac quantum algebra $A_n(R)$ over $R$, then we have :
$$\deg(\phi^{-1}) \leq \deg(\phi)^{2n - 1}$$
\end{thm}

\begin{proof} (1) If the characteristic of $R$ is positive, then the proof follows from the cited Gaber's theorem thanks to the previous second remark. 

(2) So, let's assume that the characteristic of $R$ is zero. Without loose of generality, we may also assume that $R$ is a finitely $\ZZ$-algebra. Let us assume in addition that $\deg(\phi^{-1}) > \deg(\phi)^{2n - 1}$. Let $R'$ be the the sub-algebra of the fractions field of $R$ generated by the inverses of the non zero dominating coefficients of $\phi(Y_i)$ and $\phi^{-1}(Y_i)$ for $1\leq i\leq n$, $\phi'$ the automorphism of $A_n(R')$ induced by $\phi$, $\fM'$ a maximal ideal of $R'$, $K$ the field $R'/\fM'$, and $\phi_K$ the automorphism of $A_n(K)$ induced by $\phi'$.

(3) So we have : $$\deg(\phi_K^{-1}) = \deg(\phi^{-1}) >  \deg(\phi)^{2n - 1} = \deg(\phi_K)^{2n - 1}$$

(4) On the other hand, since $R'$ is a finitely generated $\ZZ$-algebra, $K$ is finite (see for instance \cite{bourbaki}, Sect. 3, no. 4, Th. 3). So, according to (1), (3) is a contradiction.

\end{proof}

\begin{thm}[on the symplectic property of endomorphisms of Dirac quantum algebras in prime characteristic] .

If $R$ is a reduced ring of prime characteristc and if $\phi$ is an endomorphism of the Dirac quantum algebra $A_n(R)$ over $R$, then $\phi_0$ is an endomorphism of the $R$-algebra $P_n(R)$.
\end{thm}

\begin{proof}. 

(1) If $R$ is a field, then it follows from Cor. 3.3 of \cite{tsuchimoto2}. 

(2) If $R$ is a domain, then it follows from (1), considerating the the endomorphism of the $K$-algebra $K^{[2n]}$ induced by $\phi_0$, where $K$ is the fractions field of $R$.

(3) In the general case of $R$, if $p$ is the chracteristic of $R$ and $\fP$ a prime ideal of $R$, them $R/\fP$ is an $\FF_p$-algebra which is a domain. Hence the conclusion follows from (2).

\end{proof}

\begin{rem}.

(1) According to Theorem 1 above, the last sentence of Cor. 3.3 of \cite{tsuchimoto2}, according to which the the determinant of any endomorphism of a canonical Poisson algebra over a field of positive chracteristic is one, as it is well known for characteristic zero, is a mistake.

(2) But, still thanks to this Theorem 1, this mistake could be easily repeared in the proof of the main result of \cite{tsuchimoto2} which is its Cor. 7.3, proving the equivalence of complex Jacobian and Dixmier Conjectures in any dimension.
\end{rem}

\section{Reformulations in positive characteristics of Jacobian Conjecture in characteristic zero}

\begin{nota}[concerning the Classical Jacobian Conjecture in any characteristic \cite{adja1}, 3.1)].

Let us remind whith slight modifications the notations of \cite{adja2}, 3 concerning the Jacobian Conjecture in any characteristic.

\begin{enumerate}
\item
For $(n,d,p)\in \NN_*^2\times \NN$ with $p$ prime and the convention that $0$ is prime, and $K$ a field of characteristic $p$, let us denote by $CJC(n,p,d,K)$ the ``Classical Jacobian Conjecture for endomorphisms of degree at most $d$ of an algebra of polynomials in $n$ indeterminates over $K$ of characteristic $p$'' according to  \cite{adja1}, 3.1), i.e. : 

``an endomorphism of degree at most $d$ of the algebra of polynomials in $n$ indeterminates over the field $K$ of characteristic $p$ is an automorphism, if and only its jacobian is a non zero element of $K$ and induces a field extension of degree not a multiple of $p$''.
\item
For $(n,d,p)\in \NN_*^2\times\NN$ with $p$ prime, let us denote by $CJC(n,p,d)$ the ``Classical Jacobian Conjecture for endomorphisms of degree at most $d$ of an algebra of polynomials in $n$ indeterminates over a field of characteristic $p$'', i.e. the following statement : ``$CJC(n,p,d,K)$ is true for all fields $K$ of characteristic $p$''.
\item
For $(n,p)\in\NN_*\times\NN$ with $p$ prime, let us denote by $CJC(n,p)$ the ``Classical Jacobian Conjecture in $n$ determinates in characteristic $p$'', i.e. the following statement : ``$CJC(n,p,d)$ is true for all $d\in\NN_*$''.
\item
For a prime $p\in\NN$, let us denote by $CJC(p)$ the ``Classical Jacobian Conjecture in characteristic p'', i.e. the following statement : $``CJC(n,p)$ is true for all $n\in\NN_*$''.
\end{enumerate}
\end{nota}

\begin{nota}[concerning the Na\"ive Jacobian Conjecture in any characteristic].

\begin{enumerate}
\item
For $(n,p,d)\in\NN_*^2\times\NN$ with $p$ prime, and $K$ a field of characteristic $p$, let us denote by $NJC(n,p,d,K)$ the ``Na\"ive Jacobian Conjecture for endomorphisms of degree at most $d$ of an algebra of polynomials in $n$ indeterminates over $K$ of characteristic $p$'', i.e. the statement  deduced from $CJC(n,p,d,K)$ by deleting the condition on the degree of the field extension.
\item
For $(n,d,p)\in \NN_*^2\times\NN$ with $p$ prime, let us denote by $NJC(n,p,d)$ the ``Na\"ive Jacobian Conjecture  for endomorphisms of degree at most $d$ of an algebra of polynomials in $n$ indeterminates over a field of characteristic $p$'', i.e. the following statement : ``$NJC(n,p,d,K)$ is true for all fields $K$ of characteristic $p$''.
\item
For a prime $(n,p)\in\NN_*\times\NN$ with $p$ prime, let us denote by $NJC(n,p)$ the ``Na\"ive Jacobian Conjecture in $n$ determinates in characteristic $p$'', i.e. the following statement : ``$NJC(n,p,d)$ is true for all $d\in\NN_*$''.
\item
For a prime $p\in\NN$, let us denote by $NJC(p)$ the ``Na\"ive Jacobian Conjecture in characteristic p'', i.e. the following statement : $NJC(n,p)$ is true for all $n\in\NN_*$
\end{enumerate}
\end{nota}

\begin{rem}.
\begin{enumerate}
\item
Thanks to Prop. 3.7 of \cite{adja2} and to its analogous for the statement $CJC(n,p,d,K)$ where $K$ is assumed to be a commutative domain, the statements $CJC(n,p,d)$ defined above and in \cite{adja2} are equivalent, as well as the the statements $NJC(n,p,d)$ defined above and in \cite{adja2}.
\item
More generally, according to \cite{bass}, (1.1), 7, and from the Formal Invertion Theorem (see for instance \cite{bourbaki1}, Ch. 3, Sect. 4, No. 4, Prop. 5) the statements $CJC(n,p,d)$ and $NJC(n,p,d)$ defined above are equivalent to the deduced statement by replacing the assumption ``field'' by```comutative ring''.
\item
The Na\"ive Jacobian Conjecture $NJC(n,p)$ in any positive characteristic $p$ is trivially false, even for $n=1$, as proved by the classical counter-example defined by the polynomial $X - X^p$ over any field $K$ of characteristic $p$, which justify the name of ``Naive Jacobian conjecture''.
\item
But this ``na\"ivety'' could be corrected, as explained in the ``First reformulation mod p theorem for the Classical Jacobian Conjecture in characteristic 0'' of \cite{adja2}, i.e. the first reformulation in positive characteristic of this conjecture, of which we remind the statement and the proof for pedagogical considerations.
\end{enumerate}
\end{rem}

\begin{thm}[on the reformulation in positive characteristic of the Jacobian Conjecture in characteristic zero].

For any $(n,d)\in\NN_*^2$, there exists $N(n,d)\in\NN_*$ such that $CJC(n,0,d)$ is equivalent to one of the statements : ``$CJC(n,p,d)$ for all primes $p>N(n,d)$'' or ``$NJC(n,p,d)$ for all primes $p>N(n,d)$''.
\end{thm}

\begin{proof}.

According to \cite{bass}, I, (1.1)8, $CJC(n,0,d)$ is equivalent to $NJC(n,0,d,\CC)$ (``Lefschetz Principle'' for automorphisms of an algebra of polynomials over a field of characteristic zero). On the other hand, thanks to Gabber's degree bound theorem for the inverse of an automorphism of an algebra of polynomials over a field (see \cite{bass}, Cor. 1.4 of Theorem 1.5, since $NJC(n,0,d,\CC)$ is a first order proposition about the field $\CC$ and since the field $\CC$ is isomorphisc to the ultraproduct of the algebraic closures of prime finite fields accarding to the ultrafilter of the co-finite subsets of the set of non zero natural prime numbers (see for instance \cite{eklof}), it follows from \L\"os theorem that there exists an integer $N(n,d)\ge d^n$ such that $NJC(n,0,d,\CC)$ is equivalent to ``$NJC(n,p,d,\overline{\FF_p})$ for all time $p>N(n,d)$'' (see for instance \cite{eklof}, th. 3.1 and cor 3.2), and hence to ``$NJC(n,p,d,K)$ for all prime $p>N(n,d)$ and all algebraically closed field $K$ of characterstic $p$'', according to the ``elementary equivalence'' of algebraically closed fields of the same characteristic (see fo instance \cite{jensen}, ch. 1, th. 1.13). So, according to the proposition 3.4 of \cite{adja2}, claiming the equivalence of $CJC(n,p,d,K)$ and $NJC(n,p,d,K)$ for any prime $p>d^n$, $CJC(n,0,d)$ is equivalent to ``$CJC(n,p,d,K)$ for all prime $p>N(n,d)$ and all algebraically closed fields $K$ of characteristic $p$''. Finally, the conclusion follows from the proposition 3.4 of \cite{adja2} and from the following lemma.
\end{proof}

\begin{lem}.

If $K\subset L$ is a fields extension, then a $K$ linear map $f : V \rightarrow W$ between $K$-vector spaces is injective (resp. surjective) if and only if $f\otimes_K L : V\otimes_K L \rightarrow W\otimes_K L$ is injective (resp. surjective).
\end{lem}

\begin{proof}
It follows from the faithfull flatness of the free $K$-vector space $L$ (see for instance \cite{bourbaki1}, Ch. 1, Sect. 3, no. 1).
\end{proof}

\section{Reformulations in positive characteristic of Poisson Conjecture in characteristic zero}

\begin{nota}[concerning the Classical Poisson Conjecture in any characteristic].

\begin{enumerate}
\item
For $(n,d,p)\in \NN_*^2\times \NN$ with $p$ prime and the convention that $0$ is prime, and $K$ a field of characteristic $p$, let us denote by $CPC(n,p,d,K)$ the ``Classical Poisson Conjecture of index $n$ in characteristic $p$ for endomorphisms of degree at most $d$ of the $n$-th canonical Poisson $K$-algebra, i.e. : 

``Any endomorphism of degree at most $d$ of a canonical Poisson algebra of index $n$ over the field $K$ of characteristic $p\geq 0$  is an automorphism if and only if it induces a field extension of degree not a multiple of $p$ and its jacobian is a non zero element of $K$ is the case where $p\leq n$''
\item
For $(n,d,p)\in \NN_*^2\times\NN$ with $p$ prime, let us denote by $CPC(n,p,d)$ the ``Classical Poisson Conjecture of index $n$ in characteristic $p$ for endomorphisms of degree at most $d$ of the $n$-th canonical Poisson algebras over a field '', i.e. the following statement : ``$CPC(n,p,d,K)$ is true for all fields $K$ of characteristic $p$''.
\item
For $(n,p)\in\NN_*\times\NN$ with $p$ prime, let us denote by $CPC(n,p)$ the ``Classical Poisson Conjecture of index $n$ in characteristic $p$'', i.e. the following statement : ``$CPC(n,p,d)$ is true for all $d\in\NN_*$''.
\item
For a prime $p\in\NN$, let us denote by $CPC(p)$ the ``Classical Poisson Conjecture in characteristic $p$'', i.e. the following statement : ``$CPC(n,p)$ is true for all $n\in\NN_*$''.
\end{enumerate}
\end{nota}

\begin{nota}[concerning the Na\"ive Poisson Conjecture in any characteristic].

\begin{enumerate}
\item
For $(n,p,d)\in\NN_*^2\times\NN$ with $p$ prime, and $K$ a field of characteristic $p$, let us denote by $NPC(n,p,d,K)$ the ``Na\"ive Poisson Conjecture of index $n$ in characteristical $p$ for endomorphisms of degree at most $d$ of the $n$-th canonical Poisson algebra over $K$'', i.e. the statement : 

``Any endomorphism of a canonical Poisson algebra of index $n$ over the field $K$ of characteristic $p$  is an automorphism''
\item
For $(n,d,p)\in \NN_*^2\times\NN$ with $p$ prime, let us denote by $NPC(n,p,d)$ the ``Na\"ive Poisson Conjecture of index $n$ in characteristic $p$ for endomorphisms of degree at most $d$ of the $n$-th canonical Poisson algebras over a field '', i.e. the following statement : ``$NPC(n,p,d,K)$ is true for all fields $K$ of characteristic $p$''.
\item
For $(n,p)\in\NN_*\times\NN$ with $p$ prime, let us denote by $NPC(n,p)$ the ``Na\"ive Poisson Conjecture of index $n$ in characteristic $p$'', i.e. the following statement : ``$NPC(n,p,d)$ is true for all $d\in\NN_*$''.
\item
For a prime $p\in\NN$, let us denote by $NPC(p)$ the ``Na\"ive Poisson Conjecture in characteristic $p$'', i.e. the following statement : ``$CPC(n,p)$ is true for all $n\in\NN_*$''.
\end{enumerate}
\end{nota}

\begin{rem}.
\begin{enumerate}
\item
It follows from Theorem 1 above, from \cite{bass}, (1.1), 7 and from Formal Invertion Theorem (see for instance \cite{bourbaki1}, Ch. 3, Sect. 4, No. 4, Prop. 5) that the statements $CPC(n,p,d)$ and $NPC(n,p,d)$ defined above are equivalent to the deduced statement by replacing the assumption ``field'' by```comutative ring''.
\item
The Na\"ive Poisson Conjecture $NPC(n,p)$ in any positive characteristic $p$ is trivially false, even for $n=1$, as proved by the counter-example induced by the classical counter-example to the The Na\"ive Jacobian Conjecture $NJC(1,p)$ in any positive characteristic $p$, i.e. the endomorphism of $P_1(K)$ defined by the polynomials $X_1 - X_1^p$ and $X_2$ over any field $K$ of characteristic $p$, which justify the name of ``Naive Poisson conjecture''.
\item
But as for the Na\"ive Jacobian Conjecture $NJC(n,p)$,  this ``na\"ivety'' could be corrected as follows :
\end{enumerate}
\end{rem}

\begin{thm}[on the reformulation in positive characteristic of the Poisson Conjecture in characteristic zero].

For any $(n,d)\in\NN_*^2$, there exists $N(n,d)\in\NN_*$ such that $CPC(n,0,d)$ is equivalent to one of the statements : ``$CPC(n,p,d)$ for all primes $p>N(n,d)$'' or ``$NPC(n,p,d)$ for all primes $p>N(n,d)$''.
\end{thm}

\begin{proof}.

According to Gabber's degree bound theorem for the inverse of an automorphism of an algebra of polynomials over a field (see \cite{bass}, Cor. 1.4 of Theorem 1.5, $NPC(n,p,d,K)$, hence $CPC(n,0,d,K) =NPC(n,0,d,K)$, are first order proposition about the field $K$. Thanks to the the previous lemma and to the ``elementary equivalence'' of algebraically closed fields of the same characteristic (see fo instance \cite{jensen}, ch. 1, th. 1.13), it follows that $CJC(n,0,d)$ is equivalent to $NJC(n,0,d,\CC)$ (``Lefschetz Principle'' for automorphisms of a canonical Poisson algebra over a field of characteristic zero). Finally, it follows from the same arguments as in the proof of the previous theorem that there exists an integer $N(n,d)\ge d^{2n}$ satisfying the statement of the theorem to be proved.
\end{proof}

\section{Reformulations in positive characteristic of Dixmier Conjecture in characteristic zero}

\begin{nota}[concerning the Classical Dixmier Conjecture in any characteristic].

\begin{enumerate}
\item
For $(n,d,p)\in \NN_*^2\times \NN$ with $p$ prime and the convention that $0$ is prime, and $K$ a field of characteristic $p$, let us denote by $CDC(n,p,d,K)$ the ``Classical Dixmier Conjecture of index $n$ in characteristic $p$ for endomorphisms of degree at most $d$ of the $n$-th Dirac quantum $K$-algebra, i.e. : 

``An endomorphism of degree at most $d$ of a the $n$-th Dirac quantum algebra over a field $K$ of characteristic $p\geq 0$ is an automorphism if and only if its restriction to the center of this algebra induces a field extension of degree not a multiple of $p$ and the jacobian of this restriction is an non zero element of $K$ is the case where $p\leq n$''
\item
For $(n,d,p)\in \NN_*^2\times\NN$ with $p$ prime, let us denote by $CDC(n,p,d)$ the ``Classical Dixmier Conjecture of index $n$ in characteristic $p$ for endomorphisms of degree at most $d$ of the $n$-th Dirac quantum algebra over a field '', i.e. the following statement : ``$CDC(n,p,d,K)$ is true for all fields $K$ of characteristic $p$''.
\item
For $(n,p)\in\NN_*\times\NN$ with $p$ prime, let us denote by $CDC(n,p)$ the ``Classical Dixmier Conjecture of index $n$ in characteristic $p$'', i.e. the following statement : ``$CDC(n,p,d)$ is true for all $d\in\NN_*$''.
\item
For a prime $p\in\NN$, let us denote by $CDC(p)$ the ``Classical Dixmier Conjecture in characteristic $p$'', i.e. the following statement : ``$CDC(n,p)$ is true for all $n\in\NN_*$''.
\end{enumerate}
\end{nota}

\begin{nota}[concerning the Na\"ive Dixmier Conjecture in any characteristic].

\begin{enumerate}
\item
For $(n,p,d)\in\NN_*^2\times\NN$ with $p$ prime, and $K$ a field of characteristic $p$, let us denote by $NDC(n,p,d,K)$ the ``Na\"ive Dixmier Conjecture of index $n$ in characteristical $p$ for endomorphisms of degree at most $d$ of the $n$-th Dirac quantum algebra over $K$'', i.e. the statement : 

``An endomorphism of degree at most $d$ of the $n$-th Dirac quantum algebra over the field $K$ of characteristic $p$  is an automorphism''
\item
For $(n,d,p)\in \NN_*^2\times\NN$ with $p$ prime, let us denote by $NDC(n,p,d)$ the ``Na\"ive Dixmier Conjecture of index $n$ in characteristic $p$ for endomorphisms of degree at most $d$ of the $n$-th Dirac quantum algebras over a field '', i.e. the following statement : ``$NDC(n,p,d,K)$ is true for all fields $K$ of characteristic $p$''.
\item
For $(n,p)\in\NN_*\times\NN$ with $p$ prime, let us denote by $NDC(n,p)$ the ``Na\"ive Dirac Conjecture of index $n$ in characteristic $p$'', i.e. the following statement : ``$NDC(n,p,d)$ is true for all $d\in\NN_*$''.
\item
For a prime $p\in\NN$, let us denote by $NPC(p)$ the ``Na\"ive Dirac Conjecture in characteristic $p$'', i.e. the following statement : ``$CDC(n,p)$ is true for all $n\in\NN_*$''.
\end{enumerate}
\end{nota}

\begin{rem}.
\begin{enumerate}
\item
According to the jacobson property of finitely generated commutative algebras (see for instance \cite{bourbaki}, Ch. V, Sect. 3, no. 4, Th. 3), and the proof of \cite{bass}, (1.1), 7, the statements $CDC(n,p,d)$ and $NDC(n,p,d)$ defined above are equivalent to the deduced statement by replacing the assumption ``field $K$ of characteristic $p$'' by```comutative algebra $K$ over the prime field of characteristic $p$''.
\item
The Na\"ive Dirac Conjecture $NDC(n,p)$ in any positive characteristic $p$ is trivially false, even for $n=1$, as proved by the counter-example induced by the classical counter-example to the The Na\"ive Jacobian Conjecture $NJC(1,p)$ in any positive characteristic $p$, i.e. the endomorphism of $A_1(K)$ defined by its elements $Y_1 - Y_1^p$ and $Y_2$ over any field $K$ of characteristic $p$, which justifies the name of ``Naive Dixmier conjecture''.
\item
But as for the Na\"ive Jacobian Conjecture $NJC(n,p)$,  this ``na\"ivety'' could be corrected as follows :
\end{enumerate}
\end{rem}

\begin{thm}[on the reformulation in positive characteristic of Dixmier Conjecture in characteristic zero].

For any $(n,d)\in\NN_*^2$, there exists $N(n,d)\in\NN_*$ such that $CDC(n,0,d)$ is equivalent to one of the statements : ``$CDC(n,p,d)$ for all primes $p>N(n,d)$'' or ``$NDC(n,p,d)$ for all primes $p>N(n,d)$''.
\end{thm}

\begin{proof}.

According to the degree bound theorem for automorphisms of Dirac quantum algebras (see Theorem 2 above), $NDC(n,p,d,K)$, hence $CDC(n,0,d,K) =NDC(n,0,d,K)$, are first order propositions about the field $K$. Thanks to the the previous lemma ant to the ``elementary equivalence'' of algebraically closed fields of the same characteristic (see fo instance \cite{jensen}, ch. 1, th. 1.13), it follows that $CDC(n,0,d)$ is equivalent to $NDC(n,0,d,\CC)$ (``Lefschetz Principle'' for automorphisms of a Dirac quantum algebra over a field of characteristic zero). Finally, it follows from the same arguments as in the proof of theorem 4 that there exists an integer $N(n,d)\ge d^{2n}$ satisfying the statement of the theorem to be proved.
\end{proof}

\section{The equivalence of Jacobian, Poisson and Dixmier Conjectures}

\begin{thm}[the United Conjectures Theorem].

\begin{enumerate}
\item
For $(n,d,p)\in \NN_*^2\times\NN$ with $p=0$ or prime, we have the following chain of implications :
$$CJC(2n,p,d) \Rightarrow CPC(n,p,d) \Rightarrow CDC(n,p,d) \Rightarrow CJC(n,p,d)$$
\item
It follows that, for $(n,p)\in \NN_*\times\NN$ with $p=0$ or prime, we have the following chain of implications :
$$CJC(2n,p) \Rightarrow CPC(n,p) \Rightarrow CDC(n,p) \Rightarrow CJC(n,p)$$
\item
Finally, it follows that for any natural number $p=0$ or prime, we have the following chain of equivalences :
$$CJC(p) \Leftrightarrow CPC(p) \Leftrightarrow CDC(p)$$
\end{enumerate}
\end{thm}

\begin{proof}.

Let $(n,d,p)\in \NN_*^2\times\NN$ with $p=0$ or a positive prime.

(1) The implication $CJC(2n,p,d) \Rightarrow CPC(n,p,d)$ follows from Theorem 1 above.

(2) In the case where $p$ is a positive prime, the implication $CPC(n,p,d) \Rightarrow CDC(n,p,d)$ follows from Revoy Theorem 2 \cite{revoy}, Theorem 3.1 and Theorem 4.3 of \cite{ace}, and Theorem 3 above.

(3) The implication $CPC(n,0,d) \Rightarrow CDC(n,0,d)$ follows from (2) thanks to the reformulation theorems 5 and 6 above.

(2) The implication $CPC(n,p,d) \Rightarrow CDC(n,p,d)$ follows from the remark (6) of the Epilogue and to the property of extension of derivations of endomorphisms of an algebra of polynomials over a field with an invertble jacobian (see for instance \cite{wang}, Theo. 16 or \cite{matsumura}, Theo. 25.1 for a more general result).
\end{proof}

\section{Epilogue on the history of Dirac quantum algebras}

\begin{enumerate}
\item
The relations of the generators of Dirac quantum algebra $A_n(R)$ may be considered as the ``normalization'' of the famous ``commutation relations of $n$-dimensional Quantum Mechanics'' discovered in 1925 by P.A.M. (see \cite{dirac1}, Equations (11) and (12)), one year before being ``independently'' rediscovered , according to Dirac himself, (see \cite{dirac2},  second footnote of page 1) by  Born, Heisenberg and Jordan (see \cite{born}). This fundamental law of quantum mechanic postulates that, the ``moment operators'' $p_1,\ldots,p_n$ and the ``position operators'' $q_1,\ldots,q_n$ of the dynamical system of a particule of $n$ degree of freedom satisfy the relations : $$[p_i,q_j] =  \hbar \delta_{i,j}(2\pi \sqrt{-1})^{-1} Id, \, 1\leq i\leq n, \, 1\leq j\leq n$$ where $\hbar$ is the Planck constant and Id the identity operator (see also \cite{neuman}, Ch. 1, Section 2).
\item
This observation explains why Dirac, ``the incontestable father of the algebras $A_n(R)$'', i.e. the mathematician who in an historical paper \cite{dirac} published in 1926 introduced the first example $A_1(\CC)$ of such an algebra and who studied its first non trivial properties, called it ``the Quantum Algebra'' of which he proved that all derivations are ``inner ones'', i.e. of the form $ad(a) : x \mapsto [a,x], A_n(R) \rightarrow A_n(R)$, where $a$ is an element of $A_1(\CC)$. So, to pay a mirited homage to the fecundity of this historical paper and ``to give back to Cesar what belongs to Cesar'' and not to some one who published no line during all his life about the concerned algebras, it is imperative to rename what J. Dixmier called unjustly and by mistake in 1968 in \cite{dixmier} ``the Weyl algebra over $R$ of index $n$'' as ``the n-th Dirac algebra over $R$''. Otherwise, the whole mathematical community could be severely judged by chinese traditional wisdom claiming : ``who make a mistake and dont correct it makes another one''.
\item
The source of the unjustice and mistake concerning Dirac seems to be the first paper on the algebra $A_1(\CC)$ after the historical one of Dirac. It is and the paper \cite{littlewood} of Littewood publised in 1931, only five years after the one of Dirac and the part II of which is devoted essentially to other basic properties of $A_1(\CC)$ in addition of the one proved by Dirac, but never mentioning explicitely the this paper or the name of Dirac. Only some vague allusions could make a curious lecturer guess the existence of previous mathemaical papers on the subject inspired from Mathematical Physics . Indeed, the introduction of his paper begin by revealing the source of inspiration of his paper : ``In Mathematical Physics many quantities of a non-commutative nature are used. For most part these do not conform to any algebra that has been specially studied from an algebraic point of view'', omitting for the first time to mention the exception and the example of Dirac's paper. In this introduction, he continues : ``In this paper an attempt is made towards the classification of non-commutative algebras, as to include algebras with an infinite basis. A few of the simpler algebras are studied in detail, including some of the algebras used by Mathematical Physicists '', omitting for the second time to mention the pionner work of Dirac. At the very begining of part II of his part devoted to detailled study of the simpler algebras, he goes further in the implicit revelation of the source of inspiration of his paper by writing : ``In an endeavour to conform to usage in Quantum Theory, for the first algebra that we shall discuss, we call the two primitive elements $p$ and $x$. The quadratic modulus is taken to be $px - xp - 1$'', omiting the third time to mention Dirac's paper. 
\item
However, in his paper Littlewood proved interesting additional basic properties of $A_1(\CC)$, as its ``simplicity'', i.e. $A_1(\CC)$ contain no proper bilateral ideal (see \cite{littlewood}, Th. X), its ``integrity'', i.e. $A_1(\CC)$ is a domain (see \cite{littlewood}, Th. XII), its representation as the sub-algebra of an algebra of infinite matrices over $\CC$ generated by two explicit infinite matices (see \cite{littlewood}, proof of Th. XII), and the fact that all inversible elements of  $A_1(\CC)$ are complex numbers (see \cite{littlewood}, Th. XI), and explicit formulas for the commutator of two powers of the generators of $A_1(\CC)$ (see \cite{littlewood}, Th. XIII).
\item
From the point of view of mathematics and not quantum mechanics, one of the most important sources of fecundity of Dirac quantum algebras is that they can be interpretated in terms of differential operators, not only in the case of characteristic zero as it is well known (see for instance the introduction of \cite{dixmier}), but also in positive characteristic which is less known. More precisely we have the following canonical isomorphisms of $R$-algebras :
\item
$A_n(R)$ is canonically isomorphic to the ``ring of formal differential operators on $R[Y_1,\ldots,Y_n]$'', i.e. the ``ring of formal differential operators on $R[Y_1,\ldots,Y_n]$ generated by the partial derivations $\partial/\partial Y_i$ for $1\leq i\leq n$ '' in the sense of \cite{adja0}, Ch. 1, Def. 2.10, i.e. the additive group $R[Y_1,\ldots,Y_n]$'' endowed with the internal multiplication such that, for any $\alpha$ and $\beta$ in $\NN^n$ and any $a$ and $b$ in $R[Y_1,\ldots,Y_n]$, we have : $$(a Y^\alpha) (b Y^\beta) = \sum_{\lambda \in \NN^n, \lambda \leq \alpha} C_\alpha^\lambda a \partial^\lambda (b) Y^{\alpha + \beta - \lambda}$$ where for any $\lambda$ in $\NN^n$, we have : $$ Y^\lambda = Y_1^{\lambda_1}\ldots Y_n^{\lambda_n}$$ $$C_\alpha^\lambda = \prod_{1\leq i\leq n}\alpha_i !/(\lambda_i ! (\alpha_i - \lambda_i)!)$$ $$\partial^\lambda = (\partial/\partial Y_1)^{\lambda_1}\circ \ldots \circ (\partial/\partial Y_n)^{\lambda_n}$$
\item
The ``ring of formal differential operators on $R[Y_1,\ldots,Y_n]$'' in this sense is also canonically isomorphic to the ``ring of differential operators on $R[Y_1,\ldots,Y_n]$ generated by the Lie algebra of the $R[Y_1,\ldots,Y_n]$-module of derivations of $R[Y_1,\ldots,Y_n]$'' in the sense of Rinehart \cite{rinehart} (see also \cite{adja0}, Ch. 2, Exemple 1.4.4 and \cite{feldman}).
\item
From the mathematical point of view, another one of the most important sources of fecundity of Dirac quantum algebras is that they can be interpreted in terms of envelopping algebra. More precisely, J. Dixmier proved in \cite{dixmier1} in 1963 that for any bilateral ideal $I$ of the envelopping algebra $E$ of a complex nilpotent Lie algebra, $I$ is primitive if and only the complex algebra $E/I$ is isomorphic to a complex Dirac quantum algebra, increasing by this way the interest for Dirac quantum algebra, as he wrote in the introduction of \cite{dixmier}.
\item
In a following paper \cite{dixmier2} publised in 1966, Dixmier generalized the main result of the historical paper of Dirac by proving that for any positive integer $n$, all the derivations of the complex algebra $A_n(\CC)$ are inner, without any reference to the result he was generalizing.
\item
Dixmier's following paper \cite{dixmier} on the subject published in 1968 as mentioned has been until now the most influencial on on the subject, with unfortunately the greatest effect of propagation of the mistake and the unjustice concerning Dirac. This paper is devoted to detailled study of $A_1(K)$, for a field $K$ of characteristic zero, in the prolongation of Littlewood paper which is the oldest reference on the subject indicated by Dixmier in his third paper on the subject. The main result of this paper concerns the maximal commutative sub-algebras of $A_1(K)$ of which it is proved that they are finitely generated $K$-algebras of transcendence degree is 1, unless their fractions fields are necessarly pure extension of $K$. It is also given in this paper an explicit family of generators of the group of automorphisms of $A_1(K)$.
\item
Which the passing of time, one could say that the most important contribution of this last paper of Dixmier is the list a 6 problems posed at its end and the first one of which is well known now as ``Dixmier Conjecture''. The original problem 1 of this list is the following : ``Is any endomorphism of $A_1(K)$ an automorphism ?''. By extension, ``Dixmier Conjecture'' means nowadays ``any endomorphism of $A_n(K)$ is an automorphism for any interger positive $n$ and any field $K$ of characteristic zero'' (see for instance \cite{bass}, p. 297 and \cite{arno}, p. 264).
\item
Other problems of this list are more technical one and seems less fundamental problems. Problems 3 and 6 has already been solved by A. Joseph in 1975 \cite{joseph2} and Problem 5 recently by V. Bavula in 2005 \cite{bavula}, while the other 3 ones are still open when this paper was being writen.
\item
One of the most original contributions to the study of Dirac quantum algebras since the last paper of Dixmier turned to be the already cited paper \cite{revoy} of P. Revoy published in 1973 where he proved that for any field $K$ of positive characteristic and any positive integer $n$, $A_n(K)$ is an Azumaya algebra over its center $K^{[2n]}$ (see Th. 2, see also \cite{ace} for the definition of Azumaya algebras). He deduces from his result the generalization of the main theorem of Dirac historical paper to the algebra $A_n(K)$ for any positive integer $n$ and any field $K$.
\item
Among other original contributions to the study of Dirac quantum algebras, let us just mention without developpement : Gelfand-Kirillov Conjecture proposed in 1966  in \cite{gelfand} claiming that the fractions field of the envelopping algebra of a finite dimensional algebraic Lie algebra over an algebraically closed field of characteristic zero is isomorphic the the fractions field of the Dirac quantum algebra over an algebra of polynomials over this field, the proof of this conjecture for a Lie algebra of square matrices over such a field (see \cite{gelfand}, 6), for a finite dimensional nilpotent algebraic Lie algebra over such a field (see \cite{gelfand}, 7), for a finite dimensional solvable algebraic Lie algebras over such a field (see for instance \cite{joseph1}), for algebraic Lie algebras over such a field of dimension at most 8 (see \cite{alev2} and \cite{alev3}), the refutation of this conjecture for algebraic Lie algebras over such a field of dimension at least 9 (see \cite{alev1}), holonomic modules over a Dirac quantum algebra over a field of characteristic zero (see \cite{bernstein1} and \cite{bernstein2}), the ``stable range'' of such an albra and the solution of Serre Conjecture for such an algebra (see \cite{stafford}), and finally the theory of non commutative determinant of square matrices over such an algebra (see for instance \cite{adja0} and \cite{adja3}).
\end{enumerate}

\begin{center}
Authors addresses :

 {\tiny
\parbox[b]{4.5cm}{
\begin{center}
Pascal Kossivi Adjamagbo\\
Universit\'e Paris 6 -- UFR 929\\
4 Place Jussieu\\
75252 PARIS CEDEX 05 FRANCE\\
{\it e--mail: adja@math.jussieu.fr}
\end{center}}\hspace{1cm}
\parbox[b]{4.5cm}{
\begin{center}
Alexei Belov-Kanel\\
Institute of Mathematics\\
Hebrew University\\
Givat Ram, Jerusalem\\
91904 ISRAEL\\
{\it e--mail: kanel@mccme.ru}
\end{center}}\hspace{1cm}
\parbox[b]{4.5cm}{
\begin{center}
Maxim Kontsevich\\
I.H.E.S.\\
35, route de Chartres\\
91440 Bures-sur-Yvette\\
FRANCE\\
{\it e--mail: maxim@ihes.fr}
\end{center}}\hspace{1cm}
\parbox[b]{4.5cm}{
\begin{center}
Arno van den Essen\\
Department of Mathematics\\
Radboud University Nijmegen\\
Toernooiveld 6525 ED Nijmegen\\
THE NETHERLANDS\\
{\it e--mail: essen@math.ru.nl}
\end{center}}}\hspace{1cm}

\end{center}


\vspace{0.3cm}

\begin{thebibliography}{10}

\bibitem{adja0}
K. Adjamagbo,
\newblock {\em Les fondements de la th\'erie des d\'eterminants sur un domaine de Ore},
\newblock {\em  Th\`ese de Doctorat d'Etat es Sciences Math\'ematiques, Universit\'e Paris 6, Juillet 1991}.

\bibitem{adja1}
K. Adjamagbo,
\newblock {\em On separable algebras over a U.F.D. and the jacobian
conjecture in any characteristic},
\newblock {\em  in Automorphisms of Affine Spaces, Ed. van den Essen,
89-103, 1995, Kluwer Academic Publishers, Proceedings of the conference held in Cura\c cao, July 4-8 1994}.

\bibitem{adja2}
K. Adjamagbo,
\newblock {\em On isomorphisms of factorial domains and the Jacobian Conjecture in any characteristic},
\newblock Prepublication 91 (1996), Institut de Math\'ematiques de Jussieu.

\bibitem{adja3}
K. Adjamagbo,
\newblock {\em Adjamagbo Determinant and Serre Conjecture for linear groups over Dirac quantum algebras },
\newblock To appear.

\bibitem{adjawin}
K. Adjamagbo, T. Winiarski,
\newblock  {\em Refined Noether Normalization Theorem and Sharps Degree Bounds for Dominating Morphisms},
\newblock  Comm. Algebra, Vol. {\bf 33} (2005), 2387-2393.

\bibitem{ace}
K. Adjamagbo, J.-Y. Charbonnel, A. van den Essen,
\newblock {\em On ring homomorphisms of Azumaya algebras},
\newblock  see arXiv.org.

\bibitem{alev1}
J. Alev, A. Ooms, M. Van den Bergh,
\newblock  {\em A class of counterexamples to the Gelfand-Kirillov conjecture},
\newblock  Trans. Amer. Math. Soc. {\bf 348}(1982), no. 5, pp 1709--1716.

\bibitem{alev2}
J. Alev, A. Ooms, M. Van den Bergh,
\newblock  {\em The Gelfand-Kirillov conjecture for Lie algebras of dimension at most eight},
\newblock  J. Algebra {\bf 227}(2000), no. 2, pp 549--581.

\bibitem{alev3}
J. Alev, A. Ooms, M. Van den Bergh,
\newblock  {\em Corrigendum : ``The Gelfand-Kirillov conjecture for Lie algebras of dimension at most eight''},
\newblock  J. Algebra {\bf 230}(2000), no. 2, pp 749.

\bibitem{bass}
H. ~Bass, E. Connell, D. ~Wright,
\newblock  {\em The Jacobian Conjecture : Reduction of Degree and Formal Expansion of the Inverse},
\newblock  Bull. Amer. Math. Soc., Vol. {\bf 7}, Num. {\bf 2}, (1982), pp 287--330.

\bibitem{bavula}
v.v. Bavula,
\newblock {\em Dixmier's Problem 5 for the Weyl algebra},
\newblock {\em J. Algebra {\bf 283} (2005), no. 2, 604--621}.

\bibitem{bel1}
A. ~Belov, M. Kontsevich
\newblock  {\em On the stable equivalence of Jacobian and Dixmier Conjecture},
\newblock  Colloque autour des Conjectures Jacobienne et de Dixmier, Universit\'e Paris 6, 24 avril 2005, organisateurs K. Adjamagbo et J.-Y. Charbonnel.

\bibitem{bel2}
A. ~Belov, M. Kontsevich
\newblock  {\em The Jacobian Conjecture is stably equivalent to the Dixmier Conjecture},
\newblock  arXiv:math.RA/0512171 v2, 16 dec 2005.

\bibitem{bel3}
A. ~Belov, M. Kontsevich,
\newblock {\em Automorphisms of the Weyl algebra},
\newblock {\em Colloque autour des Conjectures Jacobienne et de Dixmier, Universit\'e Paris 6, 24 avril 2005, organisateurs K. Adjamagbo et J.-Y. Charbonnel, Mathematische Arbeitstagung June 10-16, 2005, Max-Planck Institut f\"ur Mathematik, Bonn, Germany,  Preprint {\bf MPIM2005-60a}, Lett. Math. Phys. {\bf 74}, no. 2, 181--199}.

\bibitem{bernstein1}
J. Bernstein,
\newblock  {\em Modules over a ring of differential operators : study of the fundamental solutions of equations with constant coefficients},
\newblock  {\em Funct. Anal. Appl., {\bf 5}, (1971), pp 89--101}.

\bibitem{bernstein2}
J. Bernstein,
\newblock  {\em The analytic continuation of generalized functions with respect to a parameter},
\newblock  {\em Funct. Anal. Appl., {\bf 6}, (1972), pp 273--285}.

\bibitem{born}
M. Born, W. Heisenberg and P. Jordan,
\newblock {\em Zur Quanten-mechanik II},
\newblock {\em Zeit. f. Phys., {\bf 35} (1926), 557}.

\bibitem{bourbaki1}
N. Bourbaki,
\newblock {\em Alg\`ebre Commutative, Ch. 1-4},
\newblock {Masson}, Paris, 1985.

\bibitem{bourbaki}
N. Bourbaki,
\newblock {\em Alg\`ebre Commutative, Ch. 5-7},
\newblock {Masson}, Paris, 1985.

\bibitem{coutinho}
S. C. Coutinho,
\newblock {\em The many avatars of a simple algebra},
\newblock {\em Amer. Math. Monthly {\bf 104} (1997), 593--604}.

\bibitem{dirac1}
P.A.M. Dirac,
\newblock {\em The Fundamental Equations of Quantum Mechanics},
\newblock {\em Roy. Soc. Proc., A, {\bf 109} (1925), 642--653}.

\bibitem{dirac2}
P.A.M. Dirac,
\newblock {\em Quantum Mechanics and a Preliminary Investigation of the Hydrogen Atom},
\newblock {\em Roy. Soc. Proc. A, {\bf 110} (1926), 561--579}.

\bibitem{dirac}
P.A.M. Dirac,
\newblock {\em On quantum algebra},
\newblock {\em Proc.  Cambridge Phil. Soc., {\bf 23} (1926), 412--418}.

\bibitem{dixmier1}
J. Dixmier,
\newblock {\em Sur les repr\'esentations unitaires des groupes de Lie nilpotents, II},
\newblock  Anais Acad. Bras. Ciens., t. {\bf 35} (1957), 325--388.

\bibitem{dixmier2}
J. Dixmier,
\newblock {\em Repr\'esentation irréductibles des alg\`ebres de Lie nilpotentes},
\newblock  J. Math. Pures et Appl., 9-i\`eme s\'erie, t. {\bf 45} (1966), 1--66.


\bibitem{dixmier}
J. Dixmier,
\newblock {\em Sur les alg\`ebres de Weyl},
\newblock  Bull. Soc. Math. France, {\bf 96} (1968), 209--242.

\bibitem{eklof}
P. Eklof,
\newblock {\em Ultraproducts for Algebraists}, in 
\newblock {\em Handbook of Mathematical Logic, p. 105-138, J. Barwise ed., North-Holland Publishing Company, Amsterdam, 1977}.

\bibitem{arno}
A.~van den Essen,
\newblock  {\em Polynomial Automorphisms and the Jacobian
Conjecture},
\newblock { Birkh\"auser Verlag} 2000.

\bibitem{arno2}
A.~van den Essen,
\newblock  {\em D-modules and the Jacobian
Conjecture},
\newblock Report 9108, University of Nijmegen, Toernooiveld, 6525 Ed Nijmegen, The Netherlands, 1991.

\bibitem{feldman}
G. L. Fel'dman,
\newblock {\em Global dimension of rings of differential operators},
\newblock Trans. Moscow Math. Soc., {\bf 1} (1982), 123--147.

\bibitem{gelfand}
I. M. Gelfand, A. A. Kirillov,
\newblock  {\em Sur les corps li\'es aux alg\`ebres enveloppantes des alg\`ebres de Lie},
\newblock  Publ. Math. I.H.E.S., tome  {\bf 31}, (1966), pp 5--19.

\bibitem{godbillon}
C. Godbillon,
\newblock  {\em Geom\'etrie differentielle et m\'ecanique classique},
\newblock { Hermann}, 1969.


\bibitem{grothendieck}
A.~Grothendieck,
\newblock  {\em R\'ecoltes et semailles},
\newblock { Universit\'e des Sciences et Techniques du Languedoc, Monpellier} 1985, and http://www.grothendieck-circle.org/.

\bibitem{jensen}
C. Jensen, H. Lenzing,
\newblock  {\em Model Theoric Algebra},
\newblock {Gordon and Breach Science Publishers}, 1989.

\bibitem{joseph1}
A. Joseph,
\newblock {\em Proof of the Gelfand-kirillov conjecture for solvable Lie algebras},
\newblock {\em Proc. Amer. Math. Soc.  {\bf 45} (1974), 1--10}.

\bibitem{joseph2}
A. Joseph,
\newblock {\em The Weyl algebra - semisimple and nilpotent elements},
\newblock {\em Amer. J. Math. {\bf 97} (1975), no. 3, 597--615}.

\bibitem{keller}
O.-H. Keller,
\newblock {\em Ganze Cremona-Transformationen},
\newblock {\em Monatsh. Math. Phys., {\bf 96} (1939), 299-306}.


\bibitem{littlewood}
D. E. Littlewood,
\newblock {\em 0n the classification of algebras},
\newblock {\em Proc. London. Math. Soc., {\bf 35} (1933), 200--240}.

\bibitem{matsumura}
H. Matsumura,
\newblock  {\em Commutative ring theory},
\newblock {Cambridge University Press}, 1989.

\bibitem{nagata}
M. Nagata,
\newblock {\em Some remarks on the two-dimensionam Jacobian Conjecture},
\newblock Chin. J. Math., {\bf 17} (1989), 1-7.


\bibitem{smale}
S. Smale,
\newblock {\em Mathematical problems for the next century},
\newblock Math. Intelligencer, {\bf 20} (1998), No. 2, 7--15.

\bibitem{revoy}
P. ~Revoy,
\newblock  {\em Alg\`ebres de Weyl en caract\'eristique $p$},
\newblock  C.R. Acad. Sci. Paris, S\'er. A-B {\bf 276} (1973), A 225-228.

\bibitem{rinehart}
G. S. Rinehart,
\newblock {\em Differential forms on general commutative algebras},
\newblock Tran. Amer. Math. Soc., {\bf 108} (1963), 195--222.

\bibitem{souriau}
J.-M. Souriau,
\newblock  {\em G\'eom\'etrie et relativit\'e},
\newblock {Hermann}, 1964.

\bibitem{stafford}
J. T. Stafford,
\newblock {\em Module structure of Weyl algebras},
\newblock J. London Math. Soc. (2), {\bf 18} (1978), 429--442.

\bibitem{tsuchimoto1}
Y. Tsuchimoto,
\newblock {\em Preliminaries on Dixmier Conjecture},
\newblock Mem. Fac. Sci. Kochi Univ., Ser. A Math. {\bf 24} (2003), 43--59.

\bibitem{tsuchimoto2}
Y. Tsuchimoto,
\newblock {\em Endomorphisms of Weyl algebra and p-curvatures},
\newblock Osaka J. Math. {\bf 42} (2005), No. 2, 435--452.

\bibitem{neuman}
John von Neumann,
\newblock {\em Mathematical Foundations of Quantum Machanics},
\newblock Princeton University Press, 1955.

\bibitem{wang}
S. S.-S. Wang,
\newblock {\em A Jacobian Criterion for Separability},
\newblock J. Algebra {\bf 65} (1980), 453--494.

\end{thebibliography}
\end{document}